\documentclass[journal,draftcls,onecolumn]{IEEEtran}
\usepackage{cite}

\ifCLASSINFOpdf
  \usepackage[pdftex]{graphicx}
\else
  \usepackage[dvips]{graphicx}
\fi
\usepackage[cmex10]{amsmath}
\usepackage{array}

\usepackage[tight,footnotesize]{subfigure}

\usepackage{amssymb,amsfonts}
\usepackage{bm}
\usepackage{epsfig}
\usepackage{enumerate}

\usepackage[all]{xy}

\newtheorem{theorem}{Theorem}
\newtheorem{lemma}{Lemma}
\newtheorem{remark}{Remark}

\newtheorem{proposition}{Proposition}

\newtheorem{problem}{Problem}

\renewcommand {\H}[0]{{\mathcal{H}}}

\newcommand {\M}[0]{{\mathcal{M}}}
\newcommand {\B}[0]{{\mathcal{B}}}

\newcommand {\Lp}[0]{{L^2(\Omega)}}
\newcommand {\Lpp}[1]{{L_{#1}^2(\Omega)}}
\newcommand {\Ha}[0]{{{H^\alpha}(\Omega)}}
\newcommand {\Hap}[1]{{H_{#1}^\alpha(\Omega)}}

\newcommand {\opt}[1]{{#1}^*}
\newcommand {\disc}[1]{\bar{#1}}
\newcommand {\optd}[1]{\opt{\disc{#1}}}
\newcommand {\dopt}[1]{\disc{#1}^{\dag}}
\newcommand {\dcopt}[1]{{#1}^{\dag}}

\newcommand {\optsol}{(\opt{x},\opt{u})}
\newcommand {\optdsol}{(\optd{x},\optd{u})}
\newcommand {\doptsol}{(\dopt{x},\dopt{u})}
\newcommand {\feassol}{(x,u)}
\newcommand {\dfeassol}{(\disc{x},\disc{u})}

\newcommand {\interpsol}{(\Interp{x},\Interp{u})}
\newcommand {\interpsollimit}{(\mathcal{I}_\infty{x},\mathcal{I}_\infty{u})}

\newcommand {\supTH}[0]{^\text{th}}
\newcommand {\Interp}[1]{\mathcal{I}_N{#1}}
\newcommand {\Deriv}[1]{\mathcal{D}_N{#1}}

\newcommand {\LGLgrid}[0]{\Gamma^\text{LGL}}

\begin{document}
%
\title{Optimal Control of Inhomogeneous Ensembles}
%
%
%

\author{Justin~Ruths~
        and~Jr-Shin~Li,~\IEEEmembership{Member,~IEEE,}
\thanks{Support for this work was provided by the National Science Foundation (Career Award 0747877) and the Air Force Office of Scientific Research (Young Investigator Award FA9550-10-1-0146).}%
\thanks{J. Ruths is with the Department of Electrical and Systems Engineering, Washington University in Saint Louis, MO, 63112 USA (e-mail: jruths@wustl.edu).}
\thanks{J.-S. Li is with the Department of Electrical and Systems Engineering, Washington University in Saint Louis, MO, 63112 USA (e-mail: jsli@seas.wustl.edu).}}
\maketitle

\begin{abstract}
Inhomogeneity, in its many forms, appears frequently in practical physical systems. Readily apparent in quantum systems, inhomogeneity is caused by hardware imperfections, measurement inaccuracies, and environmental variations, and subsequently limits the performance and efficiency achievable in current experiments. In this paper, we provide a systematic methodology to mathematically characterize and optimally manipulate inhomogeneous ensembles with concepts taken from ensemble control. In particular, we develop a computational method to solve practical quantum pulse design problems cast as optimal ensemble control problems, based on multidimensional pseudospectral approximations. We motivate the utility of this method by designing pulses for both standard and novel applications. We also show the convergence of the pseudospectral method for optimal ensemble control. The concepts developed here are applicable beyond quantum control, such as to neuron systems, and furthermore to systems with by parameter uncertainty, which pervade all areas of science and engineering.
\end{abstract}


%
\IEEEpeerreviewmaketitle

\section{Introduction}
\IEEEPARstart{R}{ecent} advancements in quantum research have enabled breakthroughs in biology, chemistry, physics, engineering, and medicine including better methods to understand the structure of macromolecules used in biochemical signaling and drug delivery, to facilitate the fast and efficient storage of information, and to yield higher resolution medical images for diagnosis and treatment of early stage cancer \cite{ernst_principles_1987, nielsen_information_2011, cavanagh_protein_2007}. Most, if not all, measurements and manipulations of quantum systems are achieved through the appropriate design of externally applied time-varying electromagnetic pulses, or controls \cite{ernst_principles_1987}. These pulses guide the system to produce a desired time-evolution or a specific terminal state. The design of such pulses is made significantly more difficult by inherent variations within the systems of interest. Inhomogeneity is one of the fundamental obstacles for the practical implementation and physical realization of quantum science and quantum technology. In classical systems the dispersions resulting from inhomogeneity is often compensated for by feedback control. Significant research effort has been employed in the area of quantum feedback control with several promising theoretical and practical discoveries in recent years \cite{mabuchi_feedback_2005, james_feedback_2006}. There is still a large portion of quantum systems for which state feedback is either impractical or difficult to achieve due to the short timescales and large state-space of quantum phenomena. These limitations motivate us to consider the open-loop synthesis of optimal pulses that achieve a desired goal while compensating for the inhomogeneity present in the quantum ensemble.

The behavior of a bulk quantum \emph{system} is the aggregate behavior of a large ensemble of individual quantum \emph{systems}. Although in isolation these individual systems, e.g. atoms, spins, qubits, etc., are fundamentally identical, in a physical system they are distinct due to different chemical and electromagnetic environments. This variation across the ensemble exhibits itself at the macroscopic level as variation in the values of parameters that characterize the dynamics of the bulk quantum system \cite{li_unified_2011}. For example, adjacent atoms within the same macromolecule shield the full strength of the applied external magnetic fields. Varied levels of shielding create a dispersion in the frequency of the quantum spins, which is observed as inhomogeneity in the value of the natural frequency of the bulk sample. In addition, hardware imperfection causes attenuation in the applied electromagnetic pulse over the ensemble and can be represented as variation in a scaling factor multiplying the applied pulse. Often several pulses are applied in sequence in order to achieve an intricate time-evolution of the system \cite{levitt_composite_1986}. Each pulse is designed assuming an exact (usually uniform) initial system state, however, in practice, the previous pulses only prepare the system to within a neighborhood of the assumed initial state. The additive error in such a pulse sequence can cause significant performance degradation.

Guiding the evolution of inhomogeneous ensembles is a central idea in the design and implementation of quantum experiments. As such, there is a rich literature of methods addressing this class of challenging problems. Initially these were intuitive or ad-hoc methods motivated by the symmetry of the state space \cite{levitt_composite_1986,hahn_spinechoes_1950}, which were then augmented with various heuristic and specifically designed techniques \cite{shaka_composite_1983,pauly_parameter_1991}. More recently pulse design problems have been cast as optimal control problems \cite{conolly_optimal_1986, skinner_application_2003, kahlet_improving_2004, khaneja_grape_2005, maximov_monotonic_2008}. Here we present a methodology that addresses the difficulties of the current methods and is easily generalizable to any inhomogeneous ensemble or uncertain system. The proposed method has both theoretical, such as convergence rates and computational complexity, and practical, such as ease of implementation and computation time, advantages.

In this article we describe a framework to pose robust quantum pulse design problems in the language of mathematical control theory with support from new theoretical concepts in ensemble control \cite{li_ensemble_2009,li_linear_2010,li_inhomogeneous_2006} and computational advances in multidimensional pseudospectral methods adapted for ensemble systems \cite{li_unified_2011,ruths_multidimensional_2011}. In a larger context, we provide a rigorous methodology to study and control inhomogeneous  ensembles or systems with parameter uncertainty from any area or application. In the following section we introduce the problem statement as well as our theoretical and computational tools. In Section \ref{sec:examples}, we take several examples from nuclear magnetic resonance (NMR) in liquids modeled by the bilinear Bloch equations, including broadband excitation in the presence of inhomogeneity, a sequence of broadband pulses robust to variation in the initial conditions, and systems with a time-varying frequency. We then provide empirical and theoretical justifications that the solutions computed using the pseudospectral method converge to solution of the original optimal control problem.

\section{Motivation \& Theory}\label{sec:theory}
In this section, we review the underlying concepts involved in our approach to design robust quantum pulses as well as the broader mathematical formulations necessary to characterize and solve such design problems. In what follows, we present a highly general model of quantum dynamics, which demonstrates the abundance of inhomogeneity in these problems and motivates studying the control of a family of parameterized systems. We then show how the notion of ensemble control is well suited for dealing with the inherent variation and uncertainty in practical quantum systems and formulate a new type of optimal control problem based on ensemble control.

\subsection{Quantum Dynamics \& Pulse Design}

The dynamics of a quantum system is given by the time-evolution of its density matrix. We consider here general dynamics in which the system may have interaction with the environment that leads to dissipation in the system state. Under the Markovian approximation, where the environment is modeled as an infinite thermostat which has constant state, the evolution of the density matrix can be written in Lindblad form in terms of the system Hamiltonian $H(t)$ and superoperator $L(\cdot)$ which model the unitary and nonunitary dynamics \cite{lindblad_generators_1976}, respectively,
\begin{equation*}
\frac{d}{dt}{\rho}= -i[H(t),\rho]-L(\rho), \qquad (\hbar=1).
\end{equation*}
\noindent The expression of the Hamiltonian has components corresponding to free evolution Hamiltonian, $H_f$, and the control Hamiltonians $H_i$,
\begin{equation*}
H(t)= H_f+\sum_{i=1}^{m}u_i(t)H_i ,
\end{equation*}
\noindent where $u_i(t)$ are externally applied electromagnetic pulses that can be used to manipulate, or guide, the evolution of the system state. Typical pulse design problems involve designing these pulses, or controls, to bring the final state of the density matrix $\rho(T)$ as close as possible to a target operator. This problem can be transformed, by taking the expectation values of the operators involved in the state transfer, to the vector-valued, bilinear control problem, $x\in\mathbb{R}^n$ and $u\in\mathbb{R}^m$ given by,
\begin{equation}
\label{eq:quantum_single}
\frac{d}{dt}{x}=\Big[\H_d+\sum_{i=1}^{m}u_i(t)\H_i\Big]x,
\end{equation}
\noindent where $\H_d\in\mathbb{R}^{n\times n}$ corresponds to the drift evolution representing $H_f$ and $L$, $\H_i\in\mathbb{R}^{n\times n}$ corresponds to the controlled evolution representing $H_i$, and $t\in[0,T]$ \cite{li_pseudospectral_2009}. While (\ref{eq:quantum_single}) accurately represents the classical interaction of magnetic fields, in practice the effective fields - and, therefore, the matrices representing the Hamiltonians $\H_d$ and $\H_i$ - show variation in magnitude due to different chemical environments and equipment errors. The system can no longer be described by a single equation but rather by a family of equations with variation in the parameters that characterize the motion, which motivates us to consider the dispersion in the dynamics as a continuum parameterized by the system values,
\begin{equation}
\label{eq:quantum_ensemble}
\frac{d}{dt}{x}(t,s)=\Big[\H_d(s)+\sum_{i=1}^{m}u_i(t)\H_i(s)\Big]x(t,s),
\end{equation}
\noindent where $s\in D\subset\mathbb{R}^d$ is a $d$-dimensional interval representing the $d$ parameters exhibiting variation \cite{ruths_multidimensional_2011}. In a more general formulation the matrices representing the Hamiltonians can be time-dependent, $\H_d=\H_d(t,s)$ and $\H_i=\H_i(t,s)$, as in the case of random fluctuations.

Designing a single set of controls (pulses) $u_i(t)$ that simultaneously steer an ensemble of dispersive systems in (\ref{eq:quantum_ensemble}) from an initial state to a desired final state is a fundamental problem in the control of quantum systems. Moreover, similar parameterized structures can be found across all areas of science and engineering, such as in neuroscience where a single stimulus is used to trigger a simultaneous firing of neuron oscillators with distinct oscillation frequencies \cite{Li_NOLCOS2010}. In these applications full state feedback, which is required in most current methods to compensate for system uncertainty, is impractical to obtain due to the sheer number of members (and states) within the ensemble. Averaged measurement is possible in some applications, however, this type of measurement restricts the forms of available feedback. It is, then, of particular importance to consider the corresponding open-loop control problem.

\subsection{Optimal Ensemble Control}

Systems as in (\ref{eq:quantum_ensemble}) motivate the study of a new class of inhomogeneous control systems. Ensemble control \cite{li_linear_2010} is a mathematical framework to characterize parameterized systems of the form,
\begin{equation} \label{eq:ensembledynamics}
	\frac{d}{dt}{x}(t,s) = F\big(t,s,x(t,s),u(t)\big), \qquad x(0,s)=x_0(s),
\end{equation}
\noindent where $x\in\mathbb{R}^n$, $u\in\mathbb{R}^m$, $s\in D\subset\mathbb{R}^d$, with $F$ and $x_0(s)$ smooth functions of their respective arguments. The significant challenge of this class of control problems originates from requiring the same open-loop control, $u(t)$ to guide the continuum of systems from an initial distribution, $x_0(s)$, to a desired final distribution, over the corresponding function space. Fundamental properties of these systems, such as controllability, are of particular interest - specifically addressing what types of inhomogeneities can be compensated for robustly. For example, it has been shown that the controllability of an ensemble of bilinear Bloch equations, used as a sample system in this paper (see Section \ref{sec:examples}), corresponds to the synthesis of appropriate polynomials \cite{li_ensemble_2009} and controllability conditions for an ensemble of time-varying linear systems are related to the Picard criterion of Fredholm integral equations of the first kind \cite{li_linear_2010}.

Subsequently, given dynamics and initial and final distributions, we seek methods to construct controls for such steering problems. As with any control problem, in general there may be many possible solutions that satisfy a state-to-state ensemble control problem. In addition there are often benefits, penalties, and limitations that are associated with the physical system, which can be used to rank the different solutions. Such practical considerations lead to considering an optimal control problem based on the ensemble dynamics in (\ref{eq:ensembledynamics}) which includes a cost functional (with terminal, $\varphi$, and running, $\mathcal{L}$, cost terms) to be minimized as well as possible endpoint and path constraints ($e$ and $g$, respectively),
\begin{align}
\min\ \ & \int_D \varphi(T,x(T,s))+\int_{0}^T \mathcal{L}(x(t,s),u(t))dt\ ds, \label{eq:ensemble_cost}\\
{\rm s.t.}\ \ & \frac{d}{dt}{x}(t,s) = F\big(t,s,x(t,s),u(t)\big), \nonumber \\
& e(x(0,s),x(T,s)) = 0,\nonumber \\
& g(x(t,s),u(t))\leq 0. \nonumber
\end{align}
\noindent An optimal nonlinear control problem of this form is, in general, analytically intractable. Computational methods are then required to solve such exceedingly complex optimal ensemble control problems. The idea from our previous work that constructing appropriate polynomials is a key tool in characterizing the controllability of ensemble systems of interest motivates the use of polynomials within the computational framework \cite{li_ensemble_2009}. Below we review the main ideas of the previously established pseudospectral method for optimal control to lay the foundation for our developed extension to optimal ensemble control problems. In Section \ref{sec:convergence} we complete this framework with a proof of convergence of this numerical method.

Without loss of generality, we consider a general continuous-time optimal control problem defined on the time interval $\Omega=[-1,1]$, which can be achieved by a simple affine transformation.
\begin{problem}[Continuous-Time Optimal Control] \label{prob:C}
	\begin{align}
	\min\ \ & J(x,u)=\varphi(x(1))+\int_{-1}^1 \mathcal{L}(x(t),u(t))dt, \label{eq:C_cost} \\
	{\rm s.t.}\ \ & \frac{d}{dt}{x}(t) = f(t,x(t),u(t)), \label{eq:C_dynamics} \\
	& e(x(-1),x(1)) = 0,\\
	& g(x(t),u(t))\leq 0, \label{eq:C_path}\\
	& \|u(t)\|_\infty \leq A,\quad u\in \Hap{m},\ \alpha>2
	\end{align}
	\noindent where $\varphi\in C^0$ is the terminal cost; the running cost, $\mathcal{L}\in C^\alpha$, where $C^\alpha$ is the space of continuous functions with $\alpha$ classical derivatives, and dynamics, $f\in C_n^{\alpha-1}$, where $C_n^{\alpha-1}$ is the space of $n$-vector valued $C^{\alpha-1}$ functions, with respect to the state, $x(t)\in\mathbb{R}^n$, and control, $u(t)\in\mathbb{R}^m$; $e$ and $g$ are terminal and path constraints, respectively; $\Hap{m}$ is the $m$-vector valued Sobolev space. The norm associated with the Sobolev space with $m=1$, $\Ha$, is given with respect to the $\Lp$ norm \cite{canuto_spectral_2006},
\end{problem}
$$\|h\|_{(\alpha)} = \bigg( \sum_{k=0}^\alpha \big|\big| h^{(k)} \big|\big|^2_2 \bigg)^{1/2}.$$

\subsection{Pseudospectral Method}
The pseudospectral method was originally developed to solve problems in fluid dynamics and since then has been successfully used for optimal control \cite{elnagar_pseudospectral_1995, ross_legendre_2004, fahroo_costate_2001} and applied to various areas \cite{li_pseudospectral_2009,Li_NOLCOS2010}. Pseudospectral discretization methods use expansions of orthogonal polynomials to approximate the states of the system and thereby inherit the spectral accuracy characteristic of orthogonal polynomial expansions (the $k^{\text{th}}$ coefficient of the expansion decreases faster than any inverse power of $k$) \cite{canuto_spectral_2006}. Through special properties, derivatives of these orthogonal polynomials can be expressed in terms of the polynomials themselves, making it possible to accurately approximate the differential equation that describes the dynamics with an algebraic relation imposed at a small number of discretization points. An appropriate choice of these discretization points, or nodes, facilitates the approximation of the states as well as ensuring accurate numerical integration through Gaussian quadrature.

As a collocation (or interpolation) method, the pseudospectral method uses Lagrange polynomials to approximate the states and controls of the optimal control problem,
\begin{eqnarray}
	\label{eq:Ix} &x(t) \approx I_N x(t) = \sum_{k=0}^N \bar{x}_k \ell_k(t), \\
	\label{eq:Iu} &u(t) \approx I_N u(t) = \sum_{k=0}^N \bar{u}_k \ell_k(t),
\end{eqnarray}
\noindent where $x(t_k) = I_N x(t_k) = \bar{x}_k$ and $u(t_k) = I_N u(t_k) = \bar{u}_k$ because the Lagrange polynomials have the property $\ell_k(t_i) = \delta_{ki}$, where $\delta_{ki}$ is the Kronecker delta function and $t_k$ are the interpolation nodes \cite{szego_orthogonal_1959}. Therefore, the coefficients $\bar{x}_k$ and $\bar{u}_k$ are the discretized values of the original problem and become the decision variables of the subsequent discrete problem.

Although the interpolation with Lagrange polynomials discretizes the original problem, we require a means to ensure that both the integral in the cost functional is computed accurately and the dynamics are obeyed. The integral can be approximated through Gauss quadrature; here we use Legendre polynomials as the orthogonal basis for the pseudospectral method. The Legendre-Gauss-Lobatto (LGL) quadrature approximation,
\begin{equation}\label{eq:gaussquad}
	\int_{-1}^{1} f(t) dt \approx \sum_{i=1}^{N} f(t_i) w_i, \qquad w_i = \int_{-1}^1 \ell_i(t)dt,
\end{equation}
\noindent is exact if the integrand $f\in\mathbb{P}_{2N-1}$ and the nodes $t_i\in\Gamma^{LGL}$, where $\mathbb{P}_{2N-1}$ denotes the set of polynomials of degree at most $2N-1$ and where $\Gamma^{LGL}=\{t_i:\dot{L}_N(t)|_{t_i}=0, i=1,\ldots, N-1\} \bigcup \{-1,1\}$ are the $N+1$ LGL nodes determined by the derivative of the $N^{\text{th}}$ order Legendre polynomial, $\dot{L}_N(t)$, and the interval endpoints \cite{canuto_spectral_2006}.

Using the LGL nodes, we can rewrite the Lagrange polynomials in terms of the orthogonal Legendre polynomials, which is critical to inherit the special derivative and spectral accuracy properties of the orthogonal polynomials despite using Lagrange interpolating polynomials. Given $t_k\in\Gamma^{LGL}$, we can express the Lagrange polynomials as \cite{boyd_chebyshev_2000},
\begin{equation*}\label{eq:lagleg}
	\ell_k(t) = \displaystyle\frac{1}{N(N+1)L_N(t_k)}\frac{(t^2-1) \dot{L}_N(t)}{t-t_k}.
\end{equation*}
\noindent The derivative of (\ref{eq:Ix}) at $t_j \in \Gamma^{LGL}$ is then,
\begin{equation}\label{eq:dinterpsc_j}
	\frac{d}{dt} I_N x(t_j) = \sum_{k=0}^N \bar{x}_k \dot{\ell}_k(t_j)=\sum_{k=0}^N D_{jk}\bar{x}_k \doteq (\Deriv{x})(t_j),
\end{equation}
\noindent where $D$ is the constant differentiation matrix \cite{gottlieb_theory_1984}.

We are now able to write the discretized optimal control problem using equations (\ref{eq:Ix}), (\ref{eq:Iu}), (\ref{eq:gaussquad}), and (\ref{eq:dinterpsc_j}). We transform the continuous-time problem to a constrained optimization,
\begin{align}
\min\ \ & \varphi(\bar{x}_N)+\sum_{i=0}^N \mathcal{L}(\bar{x}_i,\bar{u}_i)w^N_i, \nonumber\\
{\rm s.t.}\ \ & \sum_{k=0}^N D_{jk} \bar{x}_k =f(\bar{x}_j,\bar{u}_j), \label{eq:D1_dynamics} \\
& e(\bar{x}_0,\bar{x}_N) = 0,\nonumber\\
& g(\bar{x}_j,\bar{u}_j)\leq 0, \quad \forall\ {j\in\{0,1,\ldots,N\}} \nonumber.
\end{align}

\subsection{Multidimensional Pseudospectral Method} \label{sec:ensembleextension}
The pseudospectral method lends itself to a natural extension to consider the ensemble case, which we develop here. This is most readily apparent for a single parameter variation, i.e., $s\in[a,b]\subset\mathbb{R}$, however, is easily scaled to an arbitrary parameter dimension. In this basic case, the ensemble extension of (\ref{eq:Ix}) is
\begin{align} \label{eq:x2d}
x(t,s) &\approx I_{N\times N_s} x(t,s) =\sum_{k=0}^{N}\bar{x}_k(s)\ell_k(t) \approx \sum_{k=0}^{N}\left(\sum_{r=0}^{N_s}\bar{x}_{kr}\ell_r(s)\right)\ell_k(t).
\end{align}
\noindent The approximate derivative from (\ref{eq:dinterpsc_j}) at the LGL nodes in the respective $t$ and $s$ domains, $t_i\in\Gamma^{LGL}$ and $s_j\in\Gamma^{LGL}_{N_s}$, is
\begin{align} \label{eq:dx2d}
\frac{d}{dt}I_{N\times N_s}x(t_i,s_j)&=\sum_{k=0}^{N}D_{ik} \left(\sum_{r=0}^{N_s}\bar{x}_{kr}\ell_r(s_j)\right)=\sum_{k=0}^{N}D_{ik}\bar{x}_{kj},
\end{align}
where $\bar{x}_{kj}=x(t_k,s_j)$. In these equations we use a two-dimensional interpolating grid at the $N+1$ and $N_s+1$ LGL nodes in time and the parameter, respectively. For a general number of parameters, $\mathbf{s} = (s_1,s_2,\dots,s_d)^\prime\in D\subset\mathbb{R}^{\rm d},d>1$,
\begin{align} \label{eq:xdd}
x(t,\mathbf{s}) &\approx I_{N\times N_{s_1} \times \cdots \times N_{s_d}} x(t,\mathbf{s}) =\sum_{k=0}^{N}\bar{x}_k(\mathbf{s})\ell_k(t) = \sum_{k=0}^{N} \sum_{r_1=0}^{N_{s_1}} \cdots \sum_{r_d=0}^{N_{s_d}} \bar{x}_{k{r_1}\dots{r_d}} \ell_{r_d}(s_d) \cdots \ell_{r_1}(s_1) \ell_k(t).
\end{align}
and the derivative is, correspondingly, with $\mathbf{j} = (j_1,j_2,\dots,j_d)^\prime$,
\begin{align} \label{eq:dxdd}
\frac{d}{dt} I_{N\times N_{s_1} \times \cdots \times N_{s_d}} x(t,\mathbf{s}_\mathbf{j}) = \sum_{k=0}^{N}D_{ik}\bar{x}_{k{j_1}\dots{j_d}}.
\end{align}
The simplification from (\ref{eq:xdd}) to (\ref{eq:dxdd}) illustrates why the pseudospectral approximations are effective methods for ensemble control, as they mimic the lack of information in the parameter dimension. This aspect will also make the extension of the convergence proof for ensemble systems straightforward, as will be discussed in Section \ref{sec:convergence}.

\section{Examples}\label{sec:examples}

In this paper, we consider several examples based on the prototypical quantum control system described by the Bloch equations \cite{bloch_induction_1946}. The Bloch equations have been found to model a range of quantum phenomena from protein spectroscopy in nuclear magnetic resonance (NMR) \cite{ernst_principles_1987} and medical scans in magnetic resonance imaging (MRI) \cite{conolly_selective_1986} to Rabbi oscillations in quantum optics \cite{roskov_oscillations_2003}. In the following discussion, we will consider the specific application and terminology for NMR spectroscopy, however, the methods and results are easily transferred to these other areas of interest. In NMR spectroscopy, when the duration of the pulse design problem is small compared with the relaxation times ($T \ll T_1,T_2$, the characteristic longitudinal and transverse relaxation times, respectively), the evolution of spins can be well approximated as sequences of unitary rotations driven by the static magnetic field and the applied electromagnetic controls. In practice, the effective fields generating these rotations show variation across the quantum sample due to hardware imperfection and chemical shielding, which leads us to consider a range of magnetic field variations. The corresponding dimensionless Bloch equations in the rotating frame (see Appendix \ref{apdx:bloch}) are,
\begin{equation} \label{eq:bloch}
\frac{d}{dt}M(t,\omega,\epsilon)=\Big[\omega\Omega_z+\epsilon u(t)\Omega_y+\epsilon v(t)\Omega_x\Big]M(t,\omega,\epsilon),
\end{equation}
\noindent where $M(t,\omega,\epsilon)=(M_x(t,\omega,\epsilon),M_y(t,\omega,\epsilon),M_z(t,\omega,\epsilon))$ is the Cartesian magnetization vector for the parameter values $s=(\omega,\epsilon)$, $\omega\in[-B,B]\subset\mathbb{R}$, is the dispersion of natural frequencies, $\epsilon\in[1-\delta,1+\delta]$, $0<\delta<1$, is the amplitude attenuation factor, and $\Omega_\alpha\in\text{SO(3)}$ is the generator of rotation around the $\alpha$ axis. A pulse that compensates for the dispersion in frequency and is insensitive to the scaling of the applied controls is called a broadband pulse robust to inhomogeneity. In this section we consider several examples based on this model, including pulses robust not only to frequency dispersion and inhomogeneity, but also robust to uncertainty in initial conditions and time-varying frequencies.

\subsection{Robust $\pi$ Pulse}
Variation and dispersion in system dynamics pervade all physical experiments. In quantum systems, these inhomogeneities are often large enough to cause significant reduction in performance. The systematic framework we present here provides a rigorous way to frame any general pulse design problem for quantum control, as well as other areas of parameterized and uncertain systems.

A canonical problem in the control of quantum systems modeled by the Bloch equations is to design pulses that will accomplish a state-to-state transfer of the system. Such pulses, e.g., $\pi/2$ and $\pi$ pulses (accomplishing $\pi/2$ and $\pi$ rotations, respectively), are the fundamental building blocks of the pulse sequences used in many quantum experiments. Here, consider the inversion, or $\pi$, pulse that rotates the net magnetization from the equilibrium position ($+z$) to the $-z$ axis, i.e., $M(0)=(0\ 0\ 1)^\prime\rightarrow M(T)=(0\ 0\ -1)^\prime$. In the ensemble case, this goal corresponds to a uniform inversion of the spin vector across all choices of frequency and inhomogeneity. Specifically we consider the optimal ensemble control problem,
\begin{align}
\min\ \ & \int_{1-\delta}^{1+\delta} \int_{-B}^{B} M_z(T,\omega,\epsilon)\ d\omega\ d\epsilon +\int_{0}^T u^2(t)+v^2(t)\ dt, \label{eq:bloch_problem} \\
\text{s.t.}\ \ &\frac{d}{dt}M(t,\omega,\epsilon)=\Big[\omega\Omega_z+\epsilon u(t)\Omega_y+\epsilon v(t)\Omega_x\Big]M(t,\omega,\epsilon), \nonumber \\
&M(0,\omega,\epsilon) = (0\ 0\ 1)^\prime, \vphantom{\Bigl(} \nonumber\\
&\sqrt{u^2(t) + v^2(t)} \leq A,\ \forall\, t\in[0,T], \vphantom{\Bigl(} \nonumber
\end{align}
\noindent where $A$ is the maximum allowable amplitude and the cost functional serves to equally minimize the z-component of the spin vector (integrated across the ensemble) and the energy of the designed pulse.

\begin{figure}[t!]
\centering 
\includegraphics[width=1\columnwidth]{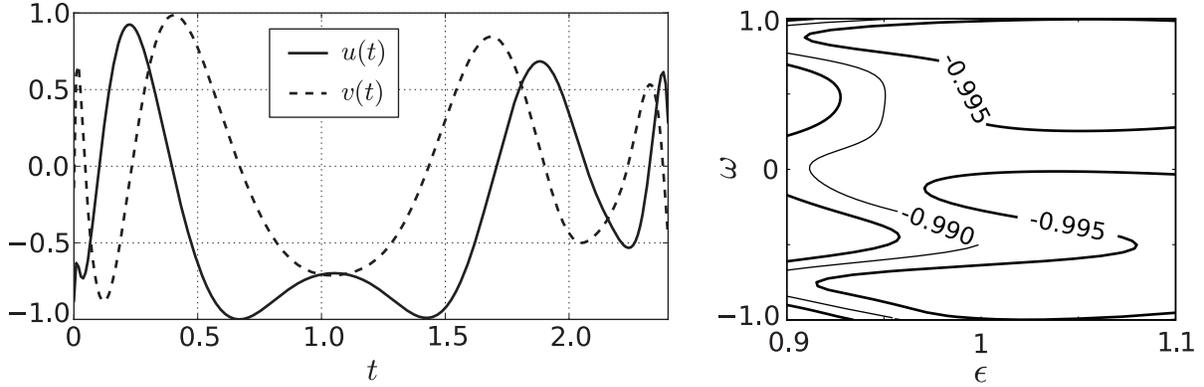}
\caption{The ``inversion'' control pulse designed by the multidimensional pseudospectral method to make the state transfer $M(0,\omega,\epsilon)=(0\ 0\ 1)^\prime\rightarrow M(T,\omega,\epsilon)=(0\ 0\ -1)^\prime$ while compensating for $\omega\in[-1,1]$ and $\epsilon\in[0.9,1.1]$, i.e., $B=1$ and $\delta=0.1$. The final states $M_z(T,\omega,\epsilon)$ shown have an average value less than -0.99, achieving a highly uniform transfer across the ensemble.}
\label{fig:pattern}
\end{figure}

Figure \ref{fig:pattern} displays a pulse that compensates for $B=1$ and $\delta=0.1$ (10\%) as well as the corresponding inversion profile. In physical units for a normalizing amplitude of 10 kHz, the maximum amplitude is $A=20$ kHz with bandwidth $\omega\in[-20,20]$ kHz and duration $T=120$ $\mu$s. Pulses developed in this manner have been implemented experimentally in true protein NMR experiments to yield significant improvement in signal recovery \cite{li_unified_2011}. Although designing individual pulses is of importance and benefit there are a myriad of other variations and uncertainties within typical quantum experiments, which calls for an approach that can address such new inhomogeneities and their corresponding challenges.

\subsection{Uncertainty in Initial Conditions}
In most experiments, individual pulses, such as the one in Figure \ref{fig:pattern}, are combined into a longer pulse sequence, which performs a more complicated manipulation of the system state with intermediate steps and goals. Even in the case of highly optimized individual pulses, as shown in the prior example, there is an error between the desired and actual final states. Moreover, pulses depend upon an exact (and usually uniform) initial condition in order to achieve their expected levels of performance. These effects combine to create a magnified accumulated error at the termination of the pulse sequence. The variation of the initial conditions of these pulses, therefore, causes significant degradation in achievable performance.

\begin{figure}[t]
\centering 
\includegraphics[width=.8\columnwidth]{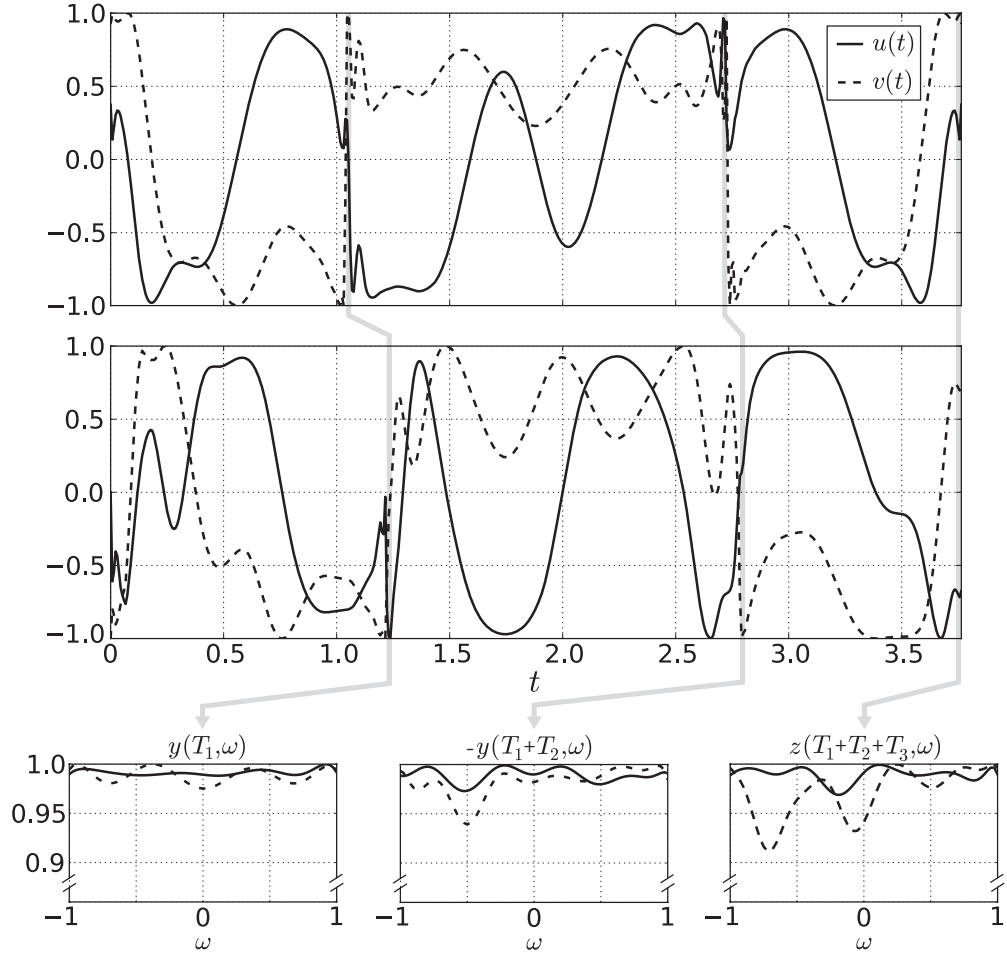}
\caption{Pulses are optimized to produce a desired ${z\rightarrow y \rightarrow -y \rightarrow z}$ evolution of the Bloch equations. The upper plot displays the concatenation of individually optimized ${z\rightarrow y}$ and ${y\rightarrow -y}$ pulses, which achieves the dashed terminal profiles shown below, with respective average performances: 0.99, 0.98, 0.97 (0.91 minimum). The middle plot displays a 3-part simultaneously-optimized pulse robust to variation in the initial condition and achieves the solid terminal profiles shown below, with respective average performances: 0.99, 0.99, 0.99 (0.97 minimum). The noticeable enhancement in performance and uniformity is due to compensating for the inhomogeneity in the initial condition of the individual pulses.}
\label{fig:initialcond}
\end{figure}

A representative example of such a pulse sequence is to perform a three step pulse sequence, which rotates the magnetization of the ensemble (1) from equilibrium ($+z$) to a point on the transverse plane (e.g. $+y$); (2) to the opposite point on the transverse plane (e.g. $-y$); (3) back to the equilibrium position ($+z$). Such pulses generally include ``phase locking'' pulses before and after the second pulse during which the magnetization dissipates. This dissipation is the portion of the experiment that is important to recover accurately and reflects a quantity to be measured, for example, a metabolic rate \cite{morris_inept_1979,payne_metabolism_2000}. If, in addition, there is accumulated error due to uncertainty in the initial conditions of the individual pulses, this leads directly to measurement inaccuracy. Here, by removing the ``phase locking'' pulses, we can abstract this pulse sequence to a unitary process and directly address any losses due to error. The controllability of the Bloch equations is shown constructing parameter-dependent (e.g. frequency, rf inhomogeneity) rotations of the spin vectors \cite{li_inhomogeneous_2006}. This, therefore, ensures that the problem with variation in initial conditions can be solved provided that the initial conditions can be parameterized by the frequency and rf-inhomogeneity.

Figure \ref{fig:initialcond} displays a three-stage optimized pulse designed by the multidimensional pseudospectral method which is robust to frequency dispersion and variation in the initial conditions of the three stages. This pulse was run as three concurrent optimizations, with the final states of one pulse fed in as the initial conditions of the next. This optimized pulse is compared with the combination of three separately optimized pulses; these combined pulses were designed with equal total duration. The terminal profiles at each intermediate goal quickly show the evidence of accumulated error in the case of the individually optimized pulses (each individual pulse has an average performance greater than 0.98). Most importantly, the uniformity of the inversion is lost in the additive error, with dips in performance down to 0.91.

\subsection{Time-Varying Frequency}
Until now, we have considered that the dispersion and uncertainty of the system are stationary. However, addressing time-varying fluctuations in parameters is also of particular theoretical and practical importance. For example, in the formulation of quantum control problems given in (\ref{eq:quantum_ensemble}) we noted that the Hamiltonians can be time-varying, motivated by such phenomena as random telegraph noise \cite{RTN}. The first step to addressing stochastic variations in such physical systems is to demonstrate control of time-varying systems, such as given by the expectation value of the corresponding random process.

\begin{figure}[t]
\centering 
\includegraphics[width=1\columnwidth]{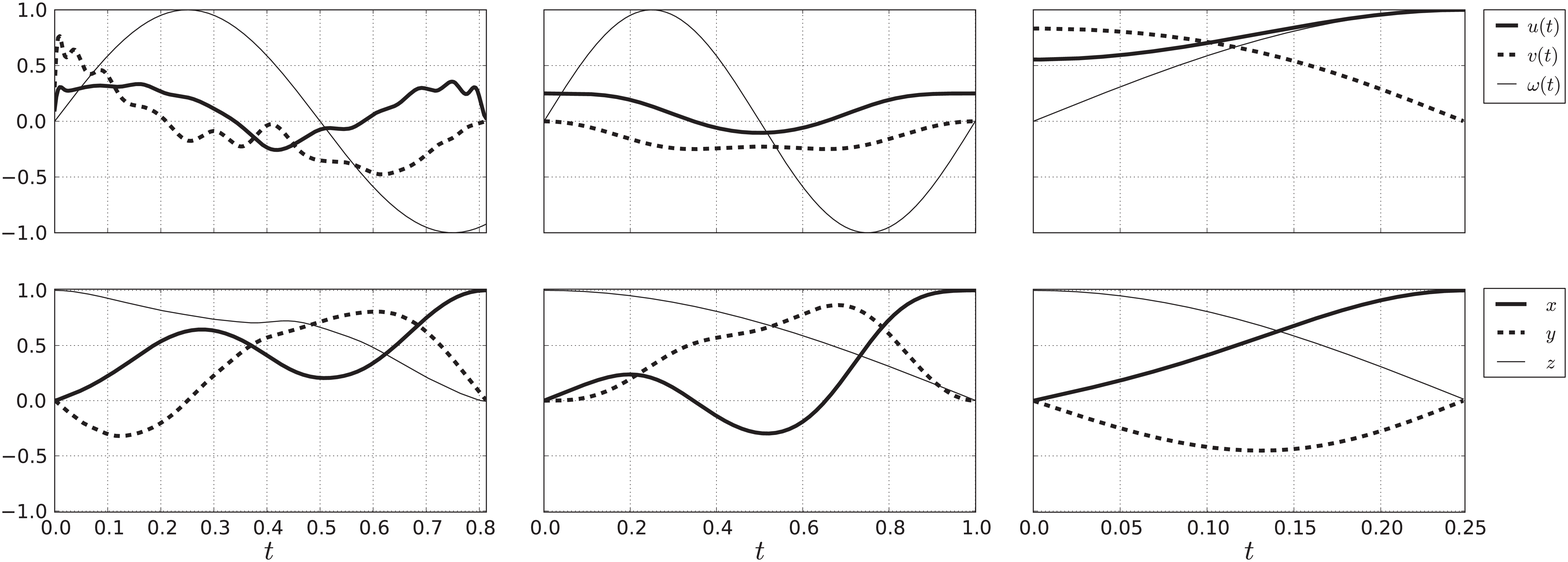}
\caption{Control pulses (top) and state trajectories (bottom) corresponding to different objectives and designed to compensate for the time-varying frequency $\omega(t)=\sin(t)$. A single-system state transfer $M(0)=(0\ 0\ 1)^\prime\rightarrow M(T)=(1\ 0\ 0)^\prime$ is designed using the terminal cost $\varphi(T) = M_x(T)$ and running costs $\mathcal{L}(t)=0$ (left), $\mathcal{L}(t)=0.1(u(t)^2+v(t)^2)$ (middle), $\mathcal{L}(t)=0.1$ (right). The terminal time was free in all cases, bounded by $T_\text{max}=1$.}
\label{fig:timevarying}
\end{figure}

Figure \ref{fig:timevarying} presents a series of optimizations designing $\pi/2$ pulses providing a state transfer $+z$ to $+x$, while compensating for a time-varying frequency, $\omega(t)=\sin(t)$. Various choices of cost functional yield different results. The arbitrary control pulse profile corresponding to the terminal cost $\varphi(T)=M_x(T)$ (Fig. \ref{fig:timevarying}, left) motivates studying optimal control methods that provide the capacity for hybrid objectives resulting in more physically meaningful controls, e.g. minimizing energy (middle) and time (right).

\section{Convergence}\label{sec:convergence}

By accepting and implementing a numerical method, we implicitly assume that the transformations and discretization used to prepare the problem for computational work does not fundamentally alter the nature of the problem. It is then critically important to show that this assumption is justified. Here we do so by both empirical and theoretical means. More specifically, we show that as the number of discretizations in the pseudospectral method (and samples in the multidimensional pseudospectral method) increases the solution of the algebraic nonlinear programming problem converges to the solution of the original continuous-time optimal control problem. For this argument, we consider a modified nonlinear programming problem statement.

\begin{problem}[Algebraic Nonlinear Programming] \label{prob:D}
	\begin{align}
	\min\ \ & \disc{J}(\disc{x},\disc{u})=\varphi(\disc{x}_N)+\sum_{k=0}^N \mathcal{L}(\disc{x}_k,\disc{u}_k)w_k \label{eq:D_cost}\\
	{\rm s.t.}\ \ & \big|\big| f\interpsol - \Deriv{x} \big|\big|_N \leq c_d N^{1-\alpha} \label{eq:D_dynamics}\\
	& e(\disc{x}_0,\disc{x}_N) = 0\\
	& g(\disc{x}_k,\disc{u}_k) \leq 0 \label{eq:D_path}\\
	& \|u_k\| \leq A \quad \forall\ k=0,1,\ldots,N
	\end{align}
	\noindent where $c_d$ is a positive constant; we define the discrete $\Lpp{n}$ norm $\|h\|_N = \sqrt{\langle h,h\rangle_N}$, for  $h,h_1,h_2\in\Lpp{n}$, $\Omega=[-1,1]$, with,
	$$\langle h_1,h_2 \rangle_N = \sum_{k=0}^N h_1^\prime(t_k) h_2(t_k)w_k,$$
	\noindent where $^\prime$ denotes the transpose and $w_k$ is the Gauss quadrature weight from (\ref{eq:gaussquad}).
\end{problem}
\begin{remark} \label{rmk:relax_dynamics}
	The dynamics in (\ref{eq:D_dynamics}) have been relaxed from the equality in (\ref{eq:D1_dynamics}) to ensure the feasibility of the discrete problem, which is used in Proposition \ref{prop:feasibility}. It is trivial to show that in the limit, as $N\rightarrow\infty$, these two conditions coincide.
\end{remark}

We seek to address three questions related to solving the continuous-time optimal control (Problem \ref{prob:C}) by solving the pseudospectral discretized constrained optimization (Problem \ref{prob:D}). Suppose a feasible solution $\feassol$ exists to Problem \ref{prob:C}. Under what conditions:
\begin{enumerate}
	\item {\bf \em Feasibility:} For a given order of approximation, $N$, does Problem \ref{prob:D} have a feasible solution, $\dfeassol$, which are the interpolation coefficients given in (\ref{eq:Ix}) and (\ref{eq:Iu})?
	\item {\bf \em Convergence:} As $N$ increases, does the sequence of optimal solutions, $\{\doptsol\}$, to Problem \ref{prob:D} have a corresponding sequence of interpolating polynomials which converges to a feasible solution of Problem \ref{prob:C}? Namely,
		$$\lim_{N\rightarrow\infty} (\Interp{\dcopt{x}},\Interp{\dcopt{u}}) = \feassol$$
	\item {\bf \em Consistency:} As $N$ increases, does the convergent sequence of interpolating polynomials corresponding to the optimal solutions of Problem \ref{prob:D} converge to an optimal solution of Problem \ref{prob:C}? Namely,
		$$\lim_{N\rightarrow\infty} (\Interp{\dcopt{x}},\Interp{\dcopt{u}}) = \optsol$$
\end{enumerate}

\begin{remark} \label{rmk:consistency}
	It is possible that Problem \ref{prob:C} has more than one optimal solution, i.e., there is more than one solution with the same optimal cost $J\optsol=J^*$. Therefore, to show that the sequence of discrete solutions converges to \emph{an} optimal solution, we can instead show that the cost of the discrete solution, $\disc{J}$, converges to the optimal cost $J^*$.
\end{remark}

Previous work has been done in the area of convergence of the pseudospectral method and we aim to augment this literature with several key insights that make convergence results applicable to a wider class of systems and relax the conditions on which the current proofs are based. Rather complete analysis has been done for the class of nonlinear systems which can be feedback linearized, including convergence rates \cite{gong_feedback_2006}. We show below that ensemble quantum systems of interest do not fall within the class of feedback linearizable systems. Work has also included general nonlinear systems, but with the assumption that the solutions of the algebraic nonlinear programming problem have a limit point (i.e., have a convergent subsequence) \cite{gong_convergence_2008}. In the language used above, this is very close to assuming ``Convergence'', which in this presentation we relax and prove Feasibility, Convergence, and Consistency directly. Finally, we examine the convergence of the multidimensional pseudospectral method as applied to ensemble optimal control problems. In what follows we consider first the convergence of the standard pseudospectral method and then discuss the convergence of the ensemble case.

We first observe that ensemble control systems of interest are not feedback linearizable \cite{isidori_nonlinear_1995}, which motivates a need for a more general convergence proof. Consider the bilinear Bloch equations in (\ref{eq:bloch}) without variation in rf inhomogeneity (i.e., $\epsilon=1$). The ability to feedback linearize a general nonlinear system is given by the Lie algebra generated by the drift and control vector fields (the conditions on this algebra must hold for each control term individually; here we consider the case for $u$). In particular, the terms $\text{ad}_{\omega\Omega_z}^{0}\Omega_y = \Omega_y$, $\text{ad}_{\omega\Omega_z}^{1}\Omega_y = -\omega\Omega_x$, $\text{ad}_{\omega\Omega_z}^{2}\Omega_y = -\omega^2\Omega_y$, $\dots$, and,
\begin{align*}
	\text{ad}_{\omega\Omega_z}^{2k-1}\Omega_y &= (-1)^k\omega^{2k-1}\Omega_x, \\
	\text{ad}_{\omega\Omega_z}^{2k}\Omega_y &= (-1)^k\omega^{2k}\Omega_y,
\end{align*}
where $k=1,2,\dots$, and $\omega$ is any value in the interval $D\subset\mathbb{R}$. It is clear that this Lie algebra, with increasing powers of the parameter $\omega$, is never closed. Therefore, the span of the appropriate Lie brackets is not involutive, which indicates that such a system is not feedback linearizable.


\subsection{Empirical Convergence}

\begin{figure}[t!]
\centering
\includegraphics[width=0.4\columnwidth]{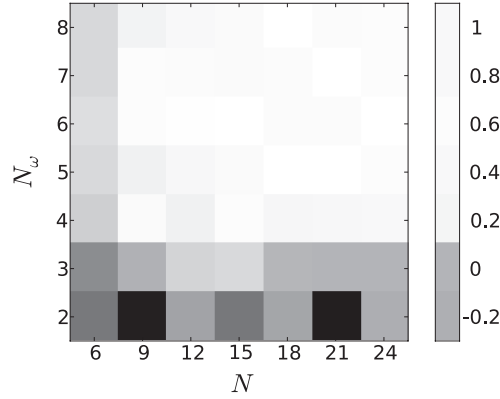}
\caption{The characteristic rapid convergence of the multidimensional pseudospectral method for a $\pi/2$ pulse designed to perform the state transfer $M(0,\omega,1)=(0\ 0\ 1)^\prime\rightarrow M(T,\omega,1)=(1\ 0\ 0)^\prime$, with $B=1$ and $T=1$. Average terminal values of $M_x(T,\omega,1)$ are shown for various choices of $N$ and $N\omega$.}
\label{fig:convergence}
\end{figure}

The orthogonal polynomials of the pseudospectral method provide spectral convergence rates similar to Fourier series approximations for periodic functions, which can easily be seen in practice. Figure \ref{fig:convergence} shows the rapid convergence of the method in both the discretization (time) and sampling (parameter) dimensions for a broadband $\pi/2$ pulse maximizing the terminal $x$ value across the ensemble. As the order of discretization (N) and/or sampling ($N_s$) increase, the method yields an objective ($\varphi(T) = M_x(T,\omega,1)$) that converges to the maximum value of unity. The low order of approximation is a characteristic of the orthogonal approximations at the heart of the numerical method. Although such empirical figures are convincing, we now show this convergence in a more rigorous fashion.

\subsection{Theoretical Preliminaries}
The results in this section will provide the foundation on which we can analyze the feasibility, convergence, and consistency of the pseudospectral approximation method for optimal control problems. We begin by presenting several key established results in polynomial approximation theory and the natural vector extensions. With these identities, we are able to then prove feasibility and convergence. We define an optimal solution to Problem \ref{prob:C} as any feasible solution that achieves the optimal cost $J\optsol=J^*$. We use this definition of an optimal solution within the subsequent preliminaries and the main result.

\begin{remark} \label{lem:xHa}
	Given Problem \ref{prob:C}, $x\in\Hap{n}$. Since $x(t)$ exists and $f\in C_n^{\alpha-1}$, all the derivatives $x^{(k)}\in C_n^0$, $\forall\ k=0,1,\dots,\alpha$ exist and are square integrable on the compact domain $\Omega$, $x^{(k)}\in\Lpp{n}$. Therefore, $x\in\Hap{n}$.
\end{remark}\vspace{1ex}

\begin{lemma}[Interpolation Error Bounds \cite{canuto_spectral_2006}, p. 289] \label{lem:poly_approx}
	If $h\in \Ha$, the following hold with $c_1,\,c_2,\,c_3,\,c>0$.
	\begin{enumerate}[(a)]
		\item \label{lem:interp} The interpolation error is bounded,
			$$\|h-\Interp{h}\|_2 \leq c_1 N^{-\alpha}\|h\|_{(\alpha)}.$$
		\item \label{lem:interp_derivative} The error between the exact derivative and the derivative of the interpolation is bounded,
			$$\|\dot{h}-\Deriv{h}\|_2 \leq c_2 N^{1-\alpha} \|h\|_{(\alpha)}.$$
			\noindent The same bound holds for the discrete $\Lp$ norm,
			$$\|\dot{h}-\Deriv{h}\|_N \leq c_3 N^{1-\alpha} \|h\|_{(\alpha)}.$$
		\item \label{lem:interp_integral} The error due to quadrature integration is bounded,
			$$\bigg| \int_{-1}^1 h(t) dt - \sum_{k=0}^N h(t_k) w_k \bigg| \leq c N^{-\alpha}\|h\|_{(\alpha)},$$
			\noindent where $t_k$ is the $k\supTH$ LGL node and $w_k$ is the corresponding $k\supTH$ weight for LGL quadrature as defined in (\ref{eq:gaussquad}).
	\end{enumerate}
\end{lemma}\vspace{1ex}

\begin{lemma} \label{lem:vector}
	If $h\in \Hap{n}$, i.e., an $n$-vector valued Sobolev space, $h = (h_1\ h_2\ \dots\ h_n)^\prime$, $h_i\in\Ha$, $i=1,2,\dots,n$.
	\begin{enumerate}[(a)]
		\item The vector-valued extension of Lemma \ref{lem:poly_approx}\ref{lem:interp} is, by the triangle inequality  on the $\Lpp{n}$ norm,
			$$\|h-\Interp{h}\|_2 \leq \sum_{i=1}^n \|h_i-\Interp{h_i}\|_2 \leq \sum_{i=1}^n c_i N^{-\alpha} \|h_i\|_{(\alpha)}.$$
		\item Similarly, \ref{lem:poly_approx}\ref{lem:interp_derivative} can be extended,
			\begin{align*}
				\|\dot{h}-\Deriv{h}\|_2 &\leq \sum_{i=1}^n \|\dot{h}-\Deriv{h}\|_2 \leq \sum_{i=1}^n c_i N^{1-\alpha} \|h_i\|_{(\alpha)} \leq c N^{1-\alpha},
			\end{align*}
			\noindent which again also holds for the discrete $\Lpp{n}$ norm.
	\end{enumerate}
\end{lemma}\vspace{1ex}

\begin{proposition}[Feasibility] \label{prop:feasibility}
	Given a solution $\feassol$ of Problem \ref{prob:C}, then Problem \ref{prob:D} has a feasible solution, $\dfeassol$, which are the corresponding interpolation coefficients.
\end{proposition}\vspace{1ex}
\begin{proof}
	Given the feasible solution $\feassol$, let $\interpsol$ be the polynomial interpolation of this solution at the LGL nodes. Our aim is to show that the coefficients of this interpolation satisfy (\ref{eq:D_dynamics})-(\ref{eq:D_path}) of Problem \ref{prob:D}. Consider the constraints imposed by the dynamics in (\ref{eq:D_dynamics}). Because the discrete norm is evaluated only at the interpolation points,
	\begin{align*}
		\|f\interpsol - \Deriv{x} \|_N &= \|f\feassol - \Deriv{x} \|_N = \|\dot{x} - \Deriv{x} \|_N \leq c_d N^{1-\alpha}
	\end{align*}
	\noindent where the last step is given by Lemma \ref{lem:vector}\ref{lem:interp_derivative}. Therefore, the interpolation coefficients $\dfeassol$ satisfy the dynamics of Problem \ref{prob:D} in (\ref{eq:D_dynamics}). We can easily show that the path constraints are also satisfied because $g(x(t),u(t))\leq 0$ for all $t\in\Omega$ by (\ref{eq:C_path}). Since this holds for all $t\in\Omega$, it also holds for all LGL nodes $t_k\in\LGLgrid$, i.e.,
	$$g(\disc{x}_k,\disc{u}_k) = g(x(t_k),u(t_k)) \leq 0,$$
	\noindent which gives (\ref{eq:D_path}). The endpoint constraints are trivially satisfied by the definition of interpolation and the presence of interpolation nodes at both endpoints. Therefore, $\dfeassol$ is a feasible solution to Problem \ref{prob:D}.
\end{proof}\vspace{1ex}

\begin{proposition}[Convergence] \label{prop:convergence}
	Given the sequence of solutions to Problem \ref{prob:D}, $\{(\disc{x},\disc{u})\}_N$, then the sequence of corresponding interpolation polynomials, $\{\interpsol\}$, has a convergent subsequence, such that
	$$\lim_{N_j\rightarrow\infty} (\mathcal{I}_{N_j}x,\mathcal{I}_{N_j}u) = \interpsollimit,$$
	\noindent which is a feasible solution to Problem \ref{prob:C}.
\end{proposition}\vspace{1ex}
\begin{proof}
	Given that $\dfeassol$ is a feasible solution of Problem \ref{prob:D}, it satisfies (\ref{eq:D_dynamics})-(\ref{eq:D_path}). Our goal is to show (i) that the sequence of solutions, $\{\interpsol\}_N$, has a convergent subsequence and (ii) that the limit point of this function subsequence is a feasible solution of Problem \ref{prob:C}, satisfying (\ref{eq:C_dynamics})-(\ref{eq:C_path}).\vspace{1ex}
	
	\noindent \emph{(i)} The sequence $\{\Interp{x}\}$, is a sequence of polynomials on a compact domain, therefore, $\Interp{x}\in\Hap{n}$. With the boundedness of the interpolating polynomials and the compactness of $\Omega$, Rellich's Theorem (see Appendix \ref{apdx:rellich}) states there is a subsequence $\{\mathcal{I}_{N_j}{x}\}$ which converges in $H_n^{\alpha-1}(\Omega)$. The same is true for the control interpolating polynomial. Therefore, there exists at least one limit point of the function sequence $\{\interpsol\}$ which we denote $\interpsollimit$.\vspace{1ex}
	
	\noindent \emph{(ii)} Explicitly writing out the calculation of the discrete norm in (\ref{eq:D_dynamics}) gives
	$$\Bigg( \sum_{k=0}^N \sum_{i=1}^n (f_i\interpsol - \Deriv{x_i})^2(t_k) \Bigg)^{1/2} \leq c_d N^{1-\alpha}.$$
	\noindent In the limit, because $f$ is continuous,
		\begin{align*}
			\lim_{N\rightarrow\infty} &\big( f_i\interpsol - \Deriv{x} \big)(t_k) = \big( f_i\interpsollimit - (\mathcal{I}_\infty{x})^\prime \big)(t_k) = 0,
			\end{align*}
	\noindent therefore,
		$$\frac{d}{dt}(\mathcal{I}_\infty{x})(t_k) = f\interpsollimit(t_k),$$
	\noindent which states that $\interpsollimit$ satisfies the dynamics in (\ref{eq:C_dynamics}) at the interpolation nodes. Moreover, as $N\rightarrow\infty$, the LGL nodes $t_k\in\LGLgrid$ are dense in $\Omega$, which further shows that $\interpsollimit$ satisfies the dynamics of Problem \ref{prob:C} at all points on the interval $\Omega$. Similarly, one can prove that this solution satisfies the path constraints because the LGL nodes become dense in $\Omega$ as $N\rightarrow\infty$ and $g(\disc{x}_k,\disc{u}_k) = g(x(t_k),u(t_k)) \leq 0$ at all LGL nodes. Again, the endpoint constraints are met exactly because the LGL grid has nodes at the endpoints.
\end{proof}\vspace{1ex}

\begin{lemma} \label{lem:interp_cost_error}
	Given $\feassol$, where $x\in\Hap{n}$, $u\in\Hap{m}$, and the corresponding interpolation coefficients, $\dfeassol$, then the error in the continuous and discrete cost functionals defined in (\ref{eq:C_cost}) and (\ref{eq:D_cost}), respectively, due to interpolation is given by,
	\begin{equation*} \label{eq:interp_cost_error} |J\feassol-\disc{J}\dfeassol| \leq cN^{-\alpha}. \end{equation*}
\end{lemma}\vspace{1ex}
\begin{remark}
	Notice that $\feassol$ and $\dfeassol$ are not required to be a feasible solutions to Problem \ref{prob:C} and \ref{prob:D}, respectively. This result characterizes the error due to interpolation.
\end{remark}\vspace{1ex}
\begin{proof} From (\ref{eq:C_dynamics}) and (\ref{eq:D_dynamics}) since $\varphi(x(1))=\varphi(\disc{x}_N)$,
	$$|J\feassol-\disc{J}\dfeassol| = \bigg| \int_{-1}^1 \mathcal{L}\feassol dt - \sum_{k=0}^N \mathcal{L}(\disc{x}_k,\disc{u}_k) w_k \bigg|.$$
	\noindent Since $\mathcal{L}\in C^\alpha$ with respect to both the state and control, $x\in\Hap{n}$ and $u\in\Hap{m}$, the composite function $\tilde{\mathcal{L}}(t)=\mathcal{L}(x(t),u(t))\in\Ha$. Let $\mathcal{L}_k=\mathcal{L}(\disc{x}_k,\disc{u}_k)$. Substituting these definitions and employing Lemma \ref{lem:poly_approx}\ref{lem:interp_integral}, we obtain
	$$\bigg| \int_{-1}^1 \tilde{\mathcal{L}}(t) dt - \sum_{k=0}^N \mathcal{L}_k w_k \bigg| \leq cN^{-\alpha} \|\tilde{\mathcal{L}}(t)\|_{(\alpha)}.$$
	\noindent Since $\tilde{\mathcal{L}}\in\Ha$, $\|\tilde{\mathcal{L}}(t)\|_{(\alpha)}$ is bounded and the result follows.
\end{proof}

\section{Main Result}
\begin{theorem}[Consistency] \label{thm:convergence}
	Suppose Problem \ref{prob:C} has an optimal solution $(\opt{x},\opt{u})$. Given a sequence of optimal solutions to Problem \ref{prob:D}, $\{(\dopt{x},\dopt{u})\}_N$, then the corresponding sequence of interpolating polynomials, $\{(\Interp{\dcopt{x}},\Interp{\dcopt{u}})\}_N$, has a limit point, $(\mathcal{I}_\infty{\dcopt{x}},\mathcal{I}_\infty{\dcopt{u}})$ which is an optimal solution to the original optimal control problem.
\end{theorem}\vspace{1ex}
\begin{proof}
	We break the proof into four sections, employing the results from the previous section.\vspace{1ex}

\noindent\emph{(i)} By Proposition \ref{prop:feasibility}, since $\optsol$ is a solution to Problem \ref{prob:C}, then for each choice of $N$, the corresponding interpolation coefficients, $\optdsol$, are a feasible solution to Problem \ref{prob:D}. By the definition of optimality of $(\dopt{x},\dopt{u})$,
		\begin{equation} \label{eq:doptimality} \disc{J}\doptsol \leq  \disc{J}\optdsol. \end{equation}
	
\noindent\emph{(ii)} By Proposition \ref{prop:convergence}, the limit point of the polynomial interpolation of the discrete optimal solution to Problem \ref{prob:D},
		$\lim_{N\rightarrow\infty} (\Interp{\dcopt{x}},\Interp{\dcopt{u}}) = (\mathcal{I}_\infty{\dcopt{x}},\mathcal{I}_\infty{\dcopt{u}})$,
	\noindent is a feasible solution of Problem \ref{prob:C}.  Therefore, we have, by the definition of the optimality of $\optsol$ and the continuity of $J$,
		\begin{align} \label{eq:coptimality}
			J\optsol &\leq  \lim_{N\rightarrow\infty} J(\Interp{\dcopt{x}},\Interp{\dcopt{u}}) = J(\mathcal{I}_\infty{\dcopt{x}},\mathcal{I}_\infty{\dcopt{u}}).
		\end{align}
	
\noindent\emph{(iii)} Using Lemma \ref{lem:interp_cost_error}, we can bound the error in the cost between the optimal solution of Problem \ref{prob:C}, $\optsol$, and the corresponding interpolating coefficients, $\optdsol$, as
		\begin{equation} \label{eq:cinterp_error} |J\optsol - \disc{J}\optdsol| \leq c_1 N^{-\alpha}. \end{equation}
		\noindent Similarly, we can bound the error in the cost between the optimal solution of Problem \ref{prob:D}, $\doptsol$, and the polynomial interpolation of this solution, $(\Interp{\dcopt{x}},\Interp{\dcopt{u}})$, as
		\begin{equation} \label{eq:dinterp_error} |J(\Interp{\dcopt{x}},\Interp{\dcopt{u}}) - \disc{J}\doptsol| \leq c_2 N^{-\alpha}. \end{equation}
		\noindent Recall that Lemma \ref{lem:interp_cost_error} does not require $(\Interp{\dcopt{x}},\Interp{\dcopt{u}})$ to be a feasible solution of Problem \ref{prob:C}. From (\ref{eq:cinterp_error}) and (\ref{eq:dinterp_error}),
		\begin{align}
			&\lim_{N\rightarrow\infty} \disc{J}\optdsol = J\optsol, \label{eq:Jlimit_copt} \\
			&\lim_{N\rightarrow\infty} \big[ J(\Interp{\dcopt{x}},\Interp{\dcopt{u}}) - \disc{J}\doptsol \big] = 0. \label{eq:Jlimit_dopt}
		\end{align}
	
\noindent\emph{(iv)} We are now ready to assemble the various pieces of this proof. Combining (\ref{eq:Jlimit_copt}) and (\ref{eq:doptimality}) we have,
	\begin{equation*} \label{eq:doptimality_limit} \lim_{N\rightarrow\infty}\disc{J}\doptsol \leq \lim_{N\rightarrow\infty}\disc{J}\optdsol = J\optsol. \end{equation*}
	\noindent Adding the result from (\ref{eq:coptimality}),
	\begin{equation} \label{eq:sandwich} \lim_{N\rightarrow\infty}\disc{J}\doptsol \leq J\optsol \leq \lim_{N\rightarrow\infty}J(\Interp{\dcopt{x}},\Interp{\dcopt{u}}). \end{equation}
	\noindent Since the difference between the left and right sides, as given by (\ref{eq:Jlimit_dopt}), decreases to zero as $N\rightarrow\infty$, the quantities $\disc{J}\doptsol$ and $J(\Interp{\dcopt{x}},\Interp{\dcopt{u}})$ converge to $J\optsol$. In particular,
	\begin{align*} \label{eq:final} 0 \leq \lim_{N\rightarrow\infty} &\big[ J\optsol - \disc{J}\doptsol\big] \leq \lim_{N\rightarrow\infty} \big[J(\Interp{\dcopt{x}},\Interp{\dcopt{u}}) - \disc{J}\doptsol\big] = 0.
	\end{align*}

\noindent Thus the optimal discrete cost $\disc{J}\doptsol$ of Problem \ref{prob:D} and the continuous cost $J(\Interp{\dcopt{x}},\Interp{\dcopt{u}})$ of the corresponding interpolation polynomials converge to the optimal cost $J\optsol$ of Problem \ref{prob:C}. Moreover, $(\mathcal{I}_\infty{\dcopt{x}},\mathcal{I}_\infty{\dcopt{u}})$ is a feasible solution to Problem \ref{prob:C} and achieves the optimal cost.  Therefore, $(\mathcal{I}_\infty{\dcopt{x}},\mathcal{I}_\infty{\dcopt{u}})$ is an optimal solution to Problem \ref{prob:C}.
\end{proof}

\begin{remark}[Ensemble Extension]\label{rmk:ensembleextension}
The nature in which the ensemble extension enters into the multidimensional pseudospectral method makes it straightforward to extend this convergence proof to the ensemble case. Section \ref{sec:ensembleextension} showed the simplicity of the derivative term in multidimensional sampling with equation (\ref{eq:dxdd}). Similarly, in the ensemble case, the constraints corresponding to the dynamics (\ref{eq:D_dynamics}) operate entirely in parallel for different parameter values. The additional integration in the cost function over the parameter domain, as in (\ref{eq:ensemble_cost}) adds another layer of quadrature approximation that can be shown to converge with arguments similar to those presented above.
\end{remark}

\section{Conclusion}
In this work we have presented a cohesive perspective and methodology for optimal control of inhomogeneous ensembles, as particularly motivated by compelling problems in quantum control and extendable to both parameterized systems in, for example, neuroscience \cite{Li_NOLCOS2010} and uncertain systems throughout science and engineering. Such systems are mathematically characterized by considering a parameterized family of differential equations indexed by a parameter vector that shows variation. Applying this rigorous framework prompts us to solve the corresponding optimal control problems with computational methods of particular form. The notion of polynomial approximation entering into the controllability analysis of the Bloch equations indicates that a modified pseudospectral method is a prime candidate. The method has natural extensions which we develop to model ensemble variation. This direct collocation method transforms the continuous-time optimal control problem into an algebraic nonlinear programming problem, which we show to be effective in a variety of applications. We supplied additional and more general arguments for the convergence of this method, in particular relaxing several assumptions and discussing the convergence characteristic of the multidimensional pseudospectral method for optimal ensemble control.

\appendices
\section{The Dimensionless Bloch Equations}\label{apdx:bloch}
The Bloch equations without relaxation, $\dot{\M}=\M\times \gamma\B_\text{eff}$, utilizes the classical description of interacting electromagnetic forces, where $\M$ is the spin magnetization vector, $\gamma$ is the gyromagnetic ratio, the effective externally applied field is $\B_\text{eff}=(B_1\cos(\omega_0 t+\phi), B_1\sin(\omega_0 t+\phi), B_0)^\prime$, $B_1(t)$ and $B_0$ are the amplitudes of the applied fields in the transverse plane and $z$ direction respectively, and $\phi(t)$ is the phase angle \cite{ernst_principles_1987}. Conventionally, the fields are given as frequencies $\gamma\B_\text{eff} = (\omega_{1x},\omega_{1y},\omega_0)$ and measured in units of Hertz. Using the generators of rotation,
\begin{equation*}
\Omega_x = \left[ \begin{array}{ccc} 0 & 0 & 0 \\ 0 & 0 & -1 \\ 0 & 1 & 0 \end{array} \right] \qquad
\Omega_y = \left[ \begin{array}{ccc} 0 & 0 & 1 \\ 0 & 0 & 0 \\ -1 & 0 & 0 \end{array} \right] \qquad
\Omega_z = \left[ \begin{array}{ccc} 0 & -1 & 0 \\ 1 & 0 & 0 \\ 0 & 0 & 0 \end{array} \right] \qquad
\end{equation*}
\noindent the Bloch equations are be given by
\begin{equation} \label{eq:blochreal}
\frac{d}{dt}\M(t)=\Big[\omega_0\Omega_z+ \omega_{1y}(t)\Omega_y+ \omega_{1x}(t)\Omega_x\Big]\M(t).
\end{equation}

\noindent If we consider variation in the applied electromagnetic fields $B_0$ and $B_1$, we can express (\ref{eq:blochreal}) in matrix form,

\begin{equation*}
	\frac{d}{dt}\left[ \begin{array}{c} \M_x(t,\omega,\epsilon) \\ \M_y(t,\omega,\epsilon) \\ \M_z(t,\omega,\epsilon) \end{array} \right] = \gamma
	\left[
		\begin{array}{ccc} 0 & -(\omega_0+\omega) & \epsilon B_1 \sin(\omega_0 t + \phi) \\
		  		\omega_0+\omega & 0 & -\epsilon B_1 \cos(\omega_0 t + \phi) \\
				-\epsilon B_1 \sin(\omega_0 t + \phi) & \epsilon B_1 \cos(\omega_0 t + \phi) & 0
		\end{array}
	\right]
	\left[ \begin{array}{c} \M_x(t,\omega,\epsilon) \\ \M_y(t,\omega,\epsilon) \\ \M_z(t,\omega,\epsilon) \end{array} \right]
\end{equation*}

\noindent where $\omega\in[-\beta,\beta]$ and $\epsilon\in[1-\delta,1+\delta]$, $0\leq\delta\leq 1$. For calculation and computation, it is useful to transform the Bloch equations into the so-called rotating frame and normalize the system by a nominal pulse amplitude $A$ to yield a dimensionless equation. Solutions based on the dimensionless equation can then be scaled for a specific choice of nominal amplitude. Consider a transformation $M = \exp(-\omega_0 \Omega_z t)\M$. In addition we scale time with $\tau=At$. It is straightforward to show that the new state equation is given by,
\begin{equation*}
\frac{d}{d\tau}M(\tau,\omega,\epsilon)=\Big[\omega\Omega_z+\epsilon u(\tau)\Omega_y+\epsilon v(\tau)\Omega_x\Big]M(\tau,\omega,\epsilon),
\end{equation*}
\noindent where $\tau\in[0,AT\times 2\pi]$, $\omega\in[-B,B]$, $B=\beta/A$, and
$$u(\tau)= \frac{\gamma B_1(\tau/A)}{A}\cos\big(\phi(\tau/A)\big) \qquad\qquad v(\tau)= \frac{\gamma B_1(\tau/A)}{A}\sin\big(\phi(\tau/A)\big),$$
\noindent (all dimensionless). Note the $2\pi$ factor in the time scaling is introduced to convert from units of Hertz to radians/second. Designing the time-varying controls $u(\tau)$ and $v(\tau)$ is equivalent to the original design of amplitude $B_1(t)$ and phase $\phi(t)$.

\section{Rellich's Theorem}\label{apdx:rellich}

\begin{theorem}[Rellich's Theorem \cite{folland_analysis}, p. 272] \label{thm:rellich}
	Suppose $\{f_k\}$ is a sequence in $H^\alpha$ such that
	\begin{enumerate}[(i)]
		\item $\sup_k\|f_k\|_{(\alpha)} < \infty$, and
		\item the $f_k$'s are all supported in a fixed compact set $V$.
	\end{enumerate}
	\noindent Then there is a subsequence $\{f_{k_j}\}$ which converges in $H^\beta$ for all $\beta < \alpha$.
\end{theorem}


\ifCLASSOPTIONcaptionsoff
  \newpage
\fi

\end{document}